# Infinite sets are non-denumerable

W. Mückenheim

University of Applied Sciences, Baumgartnerstraße 16, D-86161 Augsburg, Germany
mueckenh@rz.fh-augsburg.de

**Abstract.** Cantor's famous proof of the non-denumerability of real numbers does apply to any infinite set. The set of exclusively all natural numbers does not exist. This shows that the concept of countability is not well defined. There remains no evidence for the existence of transfinite cardinal numbers.

## 1. Introduction

Two finite sets have the same cardinality if there exists a one-to-one correspondence or bijection between them. Cantor extrapolated this theorem to include infinite sets as well: If between a set $M$ and the set $\mathbb{N}$ of all natural numbers a bijection $\mathbb{N} \leftrightarrow M$ can be established then $M$ is countable or denumerable[1] and it has the same cardinality as $\mathbb{N}$, namely $\aleph_0$. In 1874 Cantor [1] showed that the set of all algebraic numbers (including the set of all rational numbers) is denumerable, but the set $\mathbb{R}$ of all real numbers is not. He took this observation as the basis of transfinite set theory and, assuming the so-called continuum-hypothesis, he concluded that $\mathbb{R}$ has the next larger cardinal number $\aleph_1 = 2^{\aleph_0}$.

## 2. Cantor's diagonalization method

Later on Cantor simplified his proof [2] by introducing his famous diagonalization method which runs as follows: Try to set up a bijection between all real numbers $r_n$ of the interval $(0,1] \subset \mathbb{R}$ and the set $\mathbb{N}$ of all natural numbers $n$. For instance, put all the real numbers $r_n$ in an arbitrary injective sequence $n \to r_n$, i.e., a list with enumerated lines. Replace the $n^{\text{th}}$ digit $a_{nn}$ of $r_n$ by $a'_{nn} \neq a_{nn}$, avoiding $0.999... = 1.000...$ and similar identities. If possible, simply add 1. Put the changed digits together. Obtain the "diagonal number" $R = 0.a'_{11}a'_{22}a'_{33}...$ which is different from any $r_n$ of the list at least at one digit. Thus, the list cannot have been complete even for this small interval of the continuum. The assumed bijection does not exist.

---

[1] The meaning of "countable set" covers "finite set" as well as "denumerable infinite set".

Cantor's diagonalization method, however, applies to the set of natural numbers as well. It will be shown that already a bijection of ℕ onto itself is impossible.

### 3. Application to natural numbers

Following Cantor's diagonalization method, we must accept an infinite list with enumerated lines and, in case of terminating real numbers, we must extend each real number by as many zeros as are required to establish the $n^{th}$ digit (unless we assume infinitely many zeros being there at all). Of course, no such list can ever be set up, but if it is introduced in the arguing, that means, *if* we accept this procedure with respect to the sequence of real numbers, *then* we must also accept it with respect to the sequence of natural numbers.

There we obtain the same result. Consider an infinite list with lines enumerated by the sequence of the natural numbers. Add 1 to the $n^{th}$ digit (now counting from the right-hand side) of the $n^{th}$ natural number. If the number does not have enough digits, extend it by as many leading zeros as are required to establish the $n^{th}$ digit (unless we assume infinitely many zeros being there at all).

**Table 1.** A list of real numbers enumerated by the sequence of natural numbers $n$ and, in the first column, the sequence of modified natural numbers $n'$. Digits to be replaced are underlined.

| $n'$ | $n$ | $r_n$ |
|---|---|---|
| 2 | $\underline{1}$ | $0.\underline{a}_{11}a_{12}a_{13}...$ |
| 12 | $\underline{\ }2$ | $0.a_{21}\underline{a}_{22}a_{23}...$ |
| 103 | $\underline{\ \ }3$ | $0.a_{31}a_{32}\underline{a}_{33}...$ |
| ... | ... | ... |

Lining-up all the replaced digits we obtain the diagonal number $N = 111...1112$ which obviously differs from any natural number $n$ of the list by at least one digit. The question remains whether $N$ is a natural number. In order to discuss it under both possible aspects we need to distinguish the notions potential infinity and actual infinity. A quantity is *potentially infinite* (PI), if it is always finite though it has not an upper bound. An example is the function $f(n) = \{1, 2, 3, ..., n\}$ where $n \in \mathbb{N}$. A quantity is *actually infinite* (AI) if it is fixed and larger than any finite quantum. An example is $\omega = \{1, 2, 3, ...\}$.

## 4. Potential infinity

Let us first consider the list from the point of view that it is PI. Every line of the list and every digit $a'_{nn}$ of the diagonal number $R$ must have a finite index $n \in \mathbb{N}$ because the digit $a_{nn}$ of $r_n$ must be determined in order to apply the diagonalization method. The same holds for the diagonal number $N$. Each of its digits has a finite index and, therefore, a finite distance from the first digit 2. So every segment 2, 12, 112, 1112, etc. of $N$ is finite. If $N$ itself were not a finite segment, then for some digit position $n$ of $N$ being finite, the sum $N = 2 + \sum_{\nu=2}^{n} 10^{\nu}$ had to be infinite. This is impossible. So $N$ is finite and, therefore, a natural number.

## 5. Actual infinity

The view of current set theory is that every infinite set is AI. The decimal representation of any real number consists of an AI sequence of digits, and this is longer than any decimal representation of an integer. Otherwise, there was no difference between a PI sequence and an AI sequence. The AI sequence cannot consist of only finite digit positions. We see this by reflecting an AI sequence like $\pi = 3.1415...$ at the decimal point. The integer part of the reflected number ...5141.3 is not defined. As every integer with digits at only finite positions $n$ is defined, $\pi$ must have digits at infinite positions. Infinite positions, however, do not allow to determine them, let alone to find out which digits are occupying them. Therefore, these digits are not determined. In principle this would make the real numbers undetermined. Only by the factor $10^{-n}$ applied to the $n^{\text{th}}$ digit $a_{mn}$ of the real number $r_m$ and reducing its value to $a_{mn}10^{-n}$ yielding 0 in the limit, this fact becomes blurred. By the way, this property is also the reason for identities like $1 = 0.999...$ which have to be avoided in Cantor's diagonal proof. Therefore we can state that an AI list has undefined digit positions. It is insignificant by which digits these positions are occupied. The replacement of such diagonal digits $a_{nn}$ by $a'_{nn}$ is without any effect, and Cantor's diagonal argument fails.

## 6. Non-existence of ℕ

If infinity is PI, then, by diagonalization, we can always find a natural number *N* which is missing in the set of all natural numbers. If infinity is AI, then the list enumerated with all natural numbers unavoidably contains also a line enumerated by the AI number ω. Therefore the complete set ℕ of *exclusively* all natural numbers cannot and does not exist. This can already be obtained from Cantor's statement: "*Every number smaller than ω is a finite number, and its magnitude is surpassed by other finite numbers ν.*" (Jede kleinere Zahl als ω ist eine endliche Zahl und wird von anderen endlichen Zahlen ν der Größe nach übertroffen. [3]) We learn that, in addition to the set containing *every* finite number, there are *other* finite numbers.

## 7. Conclusion

Since the set ℕ of *exclusively all* natural numbers is an impossible set, we cannot construct any bijection based upon it. Any attempt to enumerate an infinite set must fail as does ℕ ↔ ℝ and already ℕ ↔ ℕ. These bijections are as impossible as $S \leftrightarrow S$, where *S* denotes the set of all sets. Therefore, infinite sets are non-denumerable, countability is not a well defined concept, transfinity is an empty expression, and with it the question concerning the cardinal number of the continuum becomes meaningless.